\documentclass[a4paper,12pt]{article}
\usepackage{makeidx}
\usepackage[T1]{fontenc}
\usepackage{amsmath,amscd,amsthm}
\usepackage{amssymb}
\usepackage{tabularx}
\usepackage{amssymb,eucal,bezier,graphicx,}
\usepackage{times,amssymb}
\usepackage[ansinew]{inputenc}
\usepackage{multicol}
\usepackage{hyperref}

\newtheorem{thm}{Theorem}[section]

\newtheorem{defi}[thm]{Definition}

\newtheorem{nota}[thm]{Note}
\newtheorem{ejem}[thm]{Example}

\begin{document}


\renewcommand{\figurename}{Figure}

\renewcommand\thefigure{\arabic{section}.\arabic{figure}} 

\numberwithin{figure}{section} 


\begin{center} {\large \bf Connecting Statistics, Probability, Algebra and Discrete Mathematics}
\end{center}

\vspace{0.10cm}

\begin{center}
{\bf Fernando López-Blázquez$^{\rm 2}$, Juan N\'u\~nez-Vald\'es$^{\rm 3}$, Silvia Recacha$^{\rm 3},$ \\ María Trinidad Villar-Li\~n\'an$^{\rm ^*, 3}$}
\end{center}

\begin{center}
{\small $^{\rm ^*}$Corresponding author,  villar@us.es} \\
{\small $^{\rm 2}$ Dpto. Estadística e Investigación Operativa. Universidad de Sevilla. \\
\small $^{\rm 3}$ Dpto. Geometr\'{\i}a y Topolog\'{\i}a. Universidad de Sevilla. \\

\small lopez@us.es, \,\, jnvaldes@us.es, \,\, silvia$\_$rg90@hotmail.com. }
\end{center}

\vspace{0.4cm}

\begin{center}
{\bf Abstract}

\end{center}

\begin{quotation}
\noindent In this paper, we connect four different branches of
Mathematics: Statistics, Probability, Algebra and Discrete
Mathematics with the objective of introducing new results on
Markov chains and evolution algebras obtained by following a
relatively new line of research, already dealt with by several
authors. It consists of the use of certain directed graphs to
facilitate the study of Markov chains and evolution algebras, as
well as to use each of the three objects to make easier the study
of the other two.
\end{quotation}

\vspace{0.1cm}

 
\vspace{0.1cm}


\vspace{0.4cm}

\section*{Introduction}
The main goal of this paper is to connect four different branches
of Mathematics to each other: Statistics and Operational Research,
Probability,  Algebra.
Indeed, we want to endow Markov chains, which are statistics
structures with two different structures more, respectively,
algebraic and discrete, apart from their already known
probabilistic character. Moreover, we also continue the research
on a relatively recent introduced topic, the {\em evolution
algebras}, to obtain new results on them, different from those
already obtained by Tian (see \cite{book, Petr}) and two of the
authors in two previous works (\cite{AMC} and \cite{NRV}), relying
on tools in principle away from these algebras as directed graphs
and Markov chains. Regarding the connection between algebraic,
discrete and statistics aspects, the reader can check
\cite{Bansaye, Staples}, for instance.

\vspace{0.1cm}

The motivation to deal with evolution algebras is the following.
At present, the study of these algebras is very booming, due to
the numerous connections between them and many other branches of
Mathematics, such as Group Theory, Dynamic Systems and Theory of Knots, among others. Furthermore,
they are also related to other sciences, including Biology, since
non-Mendelian genetics is precisely which originated them. In
fact, Tian already indicated in Chapters 2 and 4 of \cite{book}
the relationship between evolution algebras and Markov chains.

\vspace{0.1cm}

For this study we have used in this paper certain objects of
Discrete Mathematics, particularly directed graphs, thus linking
three branches of the Mathematics which at first would appear to
have no connection: evolution algebra Theory, Graph Theory and
Markov chains, so that the attainment of new properties of each of
them allows to achieve certain advances in the study of the
others.

\vspace{0.1cm}

The relationship among these three branches is based on the direct and reverse representation
of each Markov chain by an evolution algebra and the one  of the latter by
a certain graph, so that the study of the properties of each of
these objects allows us its subsequent translation to the language of the other two. This
research, which could be considered novel, was already somehow used by Tian himself
in Chapter 6 of \cite{book}.

\vspace{0.1cm}

Recently, several  works regarding the use of graphs for studying
evolution algebras have appeared in the literature
(chronologically, see \cite{AMC}, \cite{NRV}, \cite{elduque},
\cite{cabrera}, \cite{csv}, \cite{cadavid3}, \cite{cadavid2} and \cite{cadavid}, for instance). In
fact, in 2016, Tian himself invoked the researchers to continue
developing this topic in \cite{invita}.

\vspace{0.1cm}

With the results collected in this paper, for instance the novel
definitions \ref{de: graph} and \ref{de:weight} and all the
results that both involve, as Theorems \ref{th: new0}, \ref{th:new1},
\ref{th:new2}, \ref{th:new3}, \ref{th:new4} and  \ref{th:new5},
for instance, we assume the invitation made by Tian to researchers
in \cite{invita} and present some progress in the study of the
link already studied between Markov chains and evolution algebras
(see \cite{book} and \cite{cadavid}, for instance). Indeed, we
show that the properties and technics of Markov chains would allow
us to advance in the study of evolution algebras. Particularly, several
results related to evolution subalgebras are introduced and properties of generators of these algebras are analyzed.

\vspace{0.1cm}

Our intention with this contribution is to show a wide audience new mathematical models which can provide new concepts or new understanding of biological systems which the objective of finding application to multiple biological systems.

\vspace{0.1cm}

The structure of this paper is as follows: Section 1 is devoted to
recall some preliminaries on evolution algebras, Graph Theory and
Markov chains which will be used throughout the paper. In Section
2, in which the main results of the paper are shown and several examples have been incorporated to illustrate them, we deal with
the associations between pairs of these three types of objects. In
Section 3 we give the reasons for which we think that it is
indistinct to speak of evolution algebras, graphs and Markov
chains, because each property of one of them will immediately
imply the same property in the other two.

\vspace{0.1cm}

Finally, the authors would like to indicate that they have
labelled explicitly with the name of Tian in the paper all of
those definitions and results by himself, whereas the rest of
definitions and results are novel and come from our own research.

\section{Preliminaries}

In this section we recall some basic concepts related to evolution
algebras, Markov chains and Graph Theory, to be used in this
paper. We also deal with the links among all of them. For a more
general overview of these three theories, references \cite{book,
Petr} for the first, {\cite{Feller, GrimSti} for the second and
\cite{bollo, harary} for the third, among others, are available.
In this paper, we will consider finite evolution algebras and
discrete time homogeneous Markov chains with a finite number of
steps.

\subsection{Preliminaries on evolution algebras}

In this subsection, we recall some basic concepts on evolution
algebras. The most of them were introduced by Tian in \cite{book}
and other were obtained by two of the authors in \cite{AMC, NRV}.

\vspace{0.1cm}

 Evolution algebras, which  were firstly introduced
by J. P. Tian, and then jointly presented with Vojtechovsky in
2006 \cite{Petr}, and later appeared as a book by Tian in 2008
\cite{book}, are those algebras in which the relationships between
their generators $V=\{e_1, \ldots, e_n\}$ are given by
$$
\left\{
\begin{array}{l}
e_i  \cdot e_j = 0, \quad  i \neq j,  \quad 1 \leq i, j \leq n, \\
e_i^2 =e_i  \cdot e_i = \displaystyle \sum_{j = 1}^n a_{ji} \,
e_j, \quad  1 \leq i \leq n.
\end{array}
\right.
$$
where $a_{ji} \in K$ and $K$  is a field.

\vspace{0.1cm}

If $E$ is an evolution algebra over a field $K$  with a generator
set $V=\{e_i \mid i \in \Lambda\}$, the  linear map $L: E \to E$
$\mid L(e_i) = e^2_i =
 \sum_j a_{ji} e_j$, for all $i \in \Lambda$ is called the {\em evolution operator} of $E$, the coefficients $a_{ji}$ are the {\em structure constants of $E$ relative to } $V$ and the matrix $M_V:=(a_{ji})$ is said to be the {\em structure matrix of $E$ relative to $V$.}

\vspace{0.1cm}

A particular subclass of evolution algebras are the graphicable
algebras, also introduced by Tian in \cite{book}. An {\em
$n$-dimensional graphicable algebra} is a commutative, non
associative algebra, with a set of generators $V = \{e_1, e_2,
\dots, e_n \} $ endowed with relations
$$\left\{
\begin{array}{l}
e_i \cdot e_j = 0, \quad i \neq j, \quad 1 \leq i, j \leq n, \\
e_i^2 = \displaystyle \sum_{e_j \subseteq V_i}  e_j, \quad 1 \leq i \leq n. \\
\end{array}
\right.
\newline
$$

where $V_i$ is a subset of $V$.

 \noindent Thus, it is obvious that a graphicable algebra is an evolution algebra, although the converse is not true in general.

\vspace{0.1cm}

An {\em evolution subalgebra} of an evolution algebra spanned by
$\{e_1,...,e_n\}$ is a subalgebra that is spanned by $\{e_i: i \in
\Lambda \},$ for some subset $\Lambda$ of $\{1, \ldots, n\}.$

\vspace{0.1cm}

Let $E$ be an evolution algebra and $I$ be an evolution subalgebra
of $E$. It is said that $I$ is an {\em evolution ideal} of $E$ if
$E\cdot I\subseteq I$. Note that this definition implies that
every evolution subalgebra is an evolution ideal because evolution
algebras do not have an identity that characterizes them, unlike
the Lie, Malcev or Leibniz algebras, for instance.

 \vspace{0.1cm}

An evolution algebra $E$ is called {\em simple} if it has not
ideals different from $0$ and $E$, and it is called {\em
irreducible} it it has no proper subalgebra.

\vspace{0.1cm}

The following definition is given by Tian in \cite{book}: Let $E$
be an evolution algebra and $\{e_1,e_2,...e_n\}$ a set of
generators. It is said that $e_i$ {\em appears} in $x\in E$ if the
coefficient $\alpha_i \in K$ is different from $0$ in the
expression $x= \sum_{j=1}^n \alpha_j e_j$. If $e_i$ appears in
$x$, it is denoted by $e_i\prec x$.

\vspace{0.1cm}

\subsection{Preliminaries on Markov chains}
In this subsection, we recall some basic concepts on Markov chains. For further information on this topic, \cite{Karlin} and \cite{Ross} can be checked.

\vspace{0.1cm}

Given an index set $T$, a {\it stochastic process} is a collection of random variables $X=\{X(t),\, t\in T\}$. The parameter $t$ is frequently regarded as time and it may range in a discrete or a continuous set. The set of values of the variables $X(t)$ constitutes the {\it state space} of the process.

\vspace{0.1cm}

A discrete time Markov chain, $X=\{X_n,\, n=0,1,\ldots\}$ is a stochastic process satisfying the Markov property
\[
Pr[X_n=s_n\mid X_0=s_0,\ldots, X_{n-1}=s_{n-1}]=
Pr[X_n=s_n\mid  X_{n-1}=s_{n-1}]
\]
for all $n\ge 1$ and $s_i\in S$, where $S=\{s_i:\,i\in \Lambda\}$ is a numerable set of states.

\vspace{0.1cm}

A Markov chain is said to be {\it homogeneous} if the {\it transition probabilities} $p_{ij}(n)=Pr[X_{n+1}=s_j\mid X_n=s_i]$ do not depend on $n$.
In the following we will assume that our chains are homogeneous, so that the transition probabilities are denoted by $$p_{ij}=Pr[X_1=s_j\mid X_0=s_i].$$

\vspace{0.1cm}

Observe that $p_{ij}\ge 0$ for all $i,\, j\in \Lambda$ and $\sum_{j} p_{ij}=1$, for all $i\in \Lambda$. The transition probability matrix is $P=(p_{ij})_{i,j\in \Lambda}$.

\vspace{0.1cm}

A Markov chain can be regarded as a particle moving step by step in the state space. If at time $n$ the particle is at state $s_i$, then in the next step the particle will move to state $s_j$ with probability $p_{ij}$.

\vspace{0.1cm}

The $m$-step transition probability, denoted $p_{ij}^{(m)}$, is the probability that starting at $s_i$ the chain arrives at $s_j$   after exactly $m\ge 1$ steps. The Markov property implies that $P^m=(p_{ij}^{(m)})$, $m\ge 1$. By convention $p_{ij}^{(0)}=\delta_{ij}$ (Kronecker's delta).

\vspace{0.1cm}

The state $s_j$ is {\it accessible} from $s_i$ if there exists
$m\ge 0$ such that $p_{ij}^{(m)}>0$. Two states $s_i$ and $s_j$
{\it communicate} if $s_j$ is accessible from $s_i$ and $s_i$ is
accesible from $s_j$.  The relation of communication is an
equivalence relation. A Markov chain is {\it irreducible} if there
exists only one class of communicating states.

\vspace{0.1cm}

The period of a state $s_j$ is the greatest common divisor of the set of integers $m\ge 1$ such that $p_{jj}^{(m)}>0$. A state with period 1 is said to be {\it aperiodic}.

\vspace{0.1cm}

A non-empty subset of states $S_1\subset S$ is {\it closed} if $p_{ij}^{(m)}=0$, for all $i\in S_1$, $j\notin S_1$ and $m\ge 1$. If a closed subset contains only one state, then that state is {\it absorbent}.

\vspace{0.1cm}

Let $f_{ij}^{(n)}$ be the probability that the first visit to $s_j$ occurs at step $n\ge 1$ given that the chain started at $s_i$ at time $0$. The probability that a Markov chain starting at state $s_j$ returns  at some future time to $s_j$ is
$$f_{jj}=\sum_{n=1}^{\infty} f_{jj}^{(n)}.$$
If $f_{jj}=1$, then state $s_j$ is {\it recurrent}, otherwise the state is called {\it transient}.  Intuitively speaking, a state is recurrent if it is visited infinitely often. On the contrary a state is transient if after a certain moment it is not visited anymore.

\vspace{0.1cm}

From the previous results it can be deduced that the states of a Markov chain can be partitioned in classes of connected states. These classes may be closed or not. If a state $s_i$ is recurrent, then the class to which it belongs is closed. Only transient states belong to non-closed connected classes. Moreover, recurrence (or transience) is a class property, that is to say, all the states in a class of communicating states are either all recurrent, or they are all transient. Likewise, periodicity is a class property, i.e.,  all the states in the same class are either periodic with the same period, or they are all aperiodic.

\vspace{0.1cm}

Finally, note that the space of states of a Markov chain can be partitioned into two subsets, the first containing only transient states and the second only recurrent. The set of recurrent states can be further partitioned into closed connected classes.

\vspace{0.1cm}

\subsection{Preliminaries on Graph Theory}

Moving on now to Graph Theory, the most basic concepts, as graph,
directed edge, vertex, adjacency and incidence, adjacency matrix,
subgraph and walk between two vertices of a graph are assumed to
be known (\cite {harary} can be consulted for details).

\vspace{0.1cm}

A {\em weighted directed graph} $G=(V(G),E(G))$ is a pair formed by the set of vertices $V(G)$ of the graph, the set of directed edges $E(G)$ of the graph
and a function $\omega$ defined over the edges, $\omega: E(G) \to \mathbb{R}$. The image of the function in each directed edge is called the {\em weight} of the edge.

\vspace{0.1cm}

A {\em simple walk} from $u$ to $v$ in a directed graph is a sequence of vertices in which any  two consecutive vertices define a directed edge,  vertices $u$ and $v$ only appear in the beginning and at the end of the sequence, respectively, and there are no directed edges repeated.

\subsection{Links already known among evolution algebras, Graph Theory and Markov chains}

With respect to the association between evolution algebras and graphs (where {\it graph} might be of any type of them), Tian, in \cite{book}, showed how to associate a graph with an evolution algebra. He gave the following

\vspace{0.1cm}

\begin{defi}
Let $G = (V, E)$ be a directed graph, $V$ be the set of vertices and $E$ be the set of edges. It is defined the {\em associated evolution algebra with $G$} taking $V= \{e_1, e_2, \dots, e_n\}$ as the set of generators and $R$ as the set of relations of the algebra
$$
R = \left\{
\begin{array}{l}
e_i^2 = \displaystyle\sum_{e_k\in \Gamma(e_i)} e_k, \quad 1 \leq i \leq n,\\
e_i \cdot e_j =0, \quad i \neq j, \quad 1 \leq i, j \leq n.
\end{array}
\right.
\newline
$$
where $\Gamma(e_i)=\{e_k \, : \, (e_i,e_k)\in E\}$ denotes the set of vertices adjacent to $e_i$.
\end{defi}

\vspace{0.1cm}

Conversely, Tian also showed how to associate an evolution algebra
with a directed graph: he took the set of generators of the
algebra as the set of vertices and as the set of edges those
connecting the vertex $e_i$ with the vertices corresponding to
generators appearing in the expression of $e_i^2$, for each
generator $e_i$.

\vspace{0.1cm}

With respect to the relationship between graphs and Markov chains
is already known a representation of the last ones by graphs (see
\cite{Feller}, for instance) . Each homogeneous Markov chain with
the set of states $\{e_i\mid i \in \Lambda\}$ and transition
probabilities $p_{ij}=Pr[X_n=e_j\mid X_{n-1}=e_i]$ can be
associated with a weighted directed graph by taking  the set of
states as the set of vertices, the transitions of one step from a
state to other as the set of edges and the transition
probabilities as the weight of the edges. The fact that two
vertices are non adjacent means that it is not possible the
transition in one step between them.

\vspace{0.1cm}

Finally, Tian also gave in \cite{book} the relationship between
evolution algebras and Markov chains. He proved that for each
homogeneous Markov chain $X$, there is an evolution algebra $M_X$
whose structure constants are transition probabilities, and whose
generator set is the state space of the Markov chain. However, he
did not the converse procedure, that is, given an evolution algebra,
he did not define the Markov chain associated to it.

\section{Associating evolution algebras, Markov chains and graphs}

In this section we show novel results obtained on the study of the
relationships among evolution algebras, graphs and Markov chains.
All of them complete the study made by Tian in Chapter 4 of
\cite{book}.

\vspace{0.1cm}

\subsection{Some results on Markov evolution subalgebras}

In this section we show some properties and characterizations on  Markov evolution subalgebras.

\begin{defi}[Tian, Chapter 4 of \cite{book}] \label{Def21}
Let $X$ be a homogeneous Markov chain with the set of states
$\{e_i\mid i \in \Lambda\}$ and transition probabilities
$p_{ij}=Pr[X_n=e_j\mid X_{n-1}=e_i]$; The  {\em evolution algebra
$M_X$ corresponding to $X$} has $\{e_i\mid i \in
\Lambda\}$ as set of generators and the following expressions as
the laws of the algebra

$$R=
\left\{
\begin{array}{l}
e_i^2 = \displaystyle \sum_{k\in \Lambda} p_{ik}e_k, \quad i \in
\Lambda,\\
e_i \cdot e_j = 0, \quad i\neq j, \quad i, j \in
\Lambda.
\end{array}
\right.
\newline
$$
\noindent where $0\leq p_{ik}\leq 1$ and $\displaystyle
\sum_{k\in \Lambda} p_{ik}=1$.
\end{defi}

\vspace{0.1cm}

Observe that the transition probabilities matrix $(p_{ij})$ of the Markov chain $X$ defines the structure matrix of $M_X$ related
 to  $\{e_i\mid i \in \Lambda\}$.

\vspace{0.1cm}

We give now a similar definition, but in the other sense, which will allow us to set the converse result.

\begin{defi}\label{de:markovian}
Let $E$ be an evolution algebra with a set of generators $\{e_i\mid i \in \Lambda\}$ and laws given by the products

$$R=
\left\{
\begin{array}{l}
e_i^2 = \displaystyle \sum_{k\in \Lambda} a_{ik}e_k, \quad i \in
\Lambda,\\
e_i \cdot e_j = 0, \quad i\neq j, \quad i, j\in
\Lambda.
\end{array}
\right.
\newline
$$
\noindent where $0\leq a_{ik}\leq 1$ and $\sum_{k\in \Lambda}
a_{ik}=1$. The structure obtained taking the set of generators as the set of states
and the structure constants involved in the products of $E$ as
transition probabilities is called the Markov chain associated with $E.$
\end{defi}

Note that the name given to that structure is consistent, as it is proved in the following

\begin{thm} \label{th: new0}
The structure obtained in the previous definition starting from the evolution algebra $E$ is a Markov chain.
\end{thm}

\noindent {\em Proof}

\vspace{0.1cm}

According to the previous definition, it is immediate to check that the structure obtained is a Markov
chain, due to the one-to-one correspondence between $E$ and a discrete time Markov chain $X,$ with state space
the generators and transition probabilities given by the structure constants (note that each state of the Markov chain $X$ is
identified with a generator of S).   \hfill $\Box$

\begin{defi} [Tian, Introduction of Chapter 4 of \cite{book}] \label{de:markovian}
An evolution algebra which is associated with a Markov chain is
called {\em Markov evolution algebra}.
\end{defi}

\vspace{0.1cm}

Now, by construction, we will associate a weighted directed graph to a Markov evolution algebra as follows. We will call it {\em Markov graph}.

\begin{defi} \label{de: graph}
The weighted directed graph associated with a Markov evolution
algebra is $G=(V,E),$ with $V$ the set of generators of the
algebra, $(e_i, e_j) \in E$ if and only if $e_j\prec e_i^2$ and
the weight of the edge $(e_i, e_j)$ is precisely $p_{ij}$ given
by  the structure matrix of the algebra related to $V$. This
weighted directed graph is called {\em Markov graph}.
\end{defi}

\vspace{0.1cm}

Note that the Markov graph coincides with the weighted directed graph corresponding to the Markov chain which is associated with the algebra.

\begin{ejem} \label{2.4} {\rm
Let us  consider a $4$-dimensional evolution algebra defined by

$$
\begin{array}{lllllllr}
e_1 \cdot e_1 = & 0.5 \, e_1 & + & 0.2 \, e_2 &   &            & + & 0.3 \, e_4, \\
e_2 \cdot e_2 = & 0.1 \, e_1 &   &            & + & 0.9 \, e_3,&   &             \\
e_3 \cdot e_3 = &            &   &            &   & 0.4 \, e_3 & + & 0.6 \, e_4, \\
e_4 \cdot e_4 = &            &   &0.15 \, e_2 &   &            & + & 0.85\, e_4. \\
\end{array}
$$

By Definition \ref{de:markovian}, it is straightforwardly checked that  it is a Markov  evolution algebra with set of states  $\{e_1,e_2,e_3,e_4\}$, and its transition probabilities matrix is the  matrix $P$ given in  Figure \ref{image1}. Moreover, the weighted directed graph associated to such an
algebra is also shown in that figure.

}

\begin{figure}[ht!]
\centering
\includegraphics[scale=1]{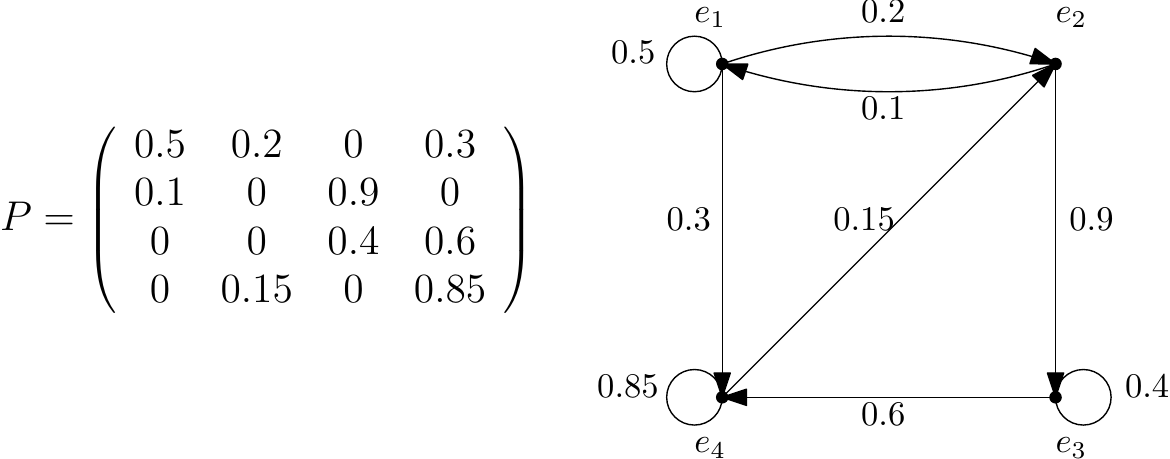}
\caption{\small The directed graph for Example \ref{2.4} and its transition probabilities matrix.} \label{image1}
\end{figure}
\end{ejem}

\vspace{0.1cm}

However, it is clear that there are  evolution algebras which are not Markov. The  next example shows one of them since it does not satisfy  Definition \ref{de:markovian}.

\vspace{0.1cm}

\begin{ejem} \label{EX} {\rm
Let $E$ be a $3$-dimensional evolution algebra defined by

$$
\begin{array}{lrrllr}
e_1 \cdot e_1 = & 0.25 \, e_1 & + & 0.3 \, e_2 & + & 0.82 \, e_3,\\
e_2 \cdot e_2 = &             & + &  0.37 \, e_2 & + & 0.63 \, e_3,\\
e_3 \cdot e_3 = & 1.3 \, e_1. &   &            &   & \\
\end{array}
$$

\vspace{0.1cm}

}

\end{ejem}

\vspace{0.1cm}

As a particular case, some graphicable algebras also have associated Markov chains. Those algebras, graphicable and Markov, are those whose laws are

$$\left\{
\begin{array}{l}
e_i^2 = e_{k(i)}, \quad \text{for  }i, k(i) \in\Lambda,\\
e_i \cdot e_j = 0, \quad i\neq j, \quad i, j \in
\Lambda.
\end{array}
\right.
\newline
$$

\vspace{0.1cm}

It can be found in \cite{AMC} the following characterization of the
evolution operator of graphicable algebras

\begin{thm}{\rm (\cite{AMC})}
Let $G$ be a simple graph with $V(G)=\{x_1,x_2,...,x_n\}$ and let $L$ be the evolution operator of graphicable algebra $A(G)$.
If we express  $L^n(e_i)=n_{i1}e_1+n_{i2}e_2+...+n_{ir}e_r$, then $n_{ij}$ coincides with the number of walks with length $n$ between the vertices corresponding to generators $e_i$ and $e_j$.
\end{thm}

The natural question is now which would be the translation of this
result in the language of Markov chains. Let us first see an
example which illustrates our intention

\vspace{0.1cm}

\begin{ejem} \label{2.5}
    {\rm
Let $X$ be an homogeneous Markov chain with set of states $\{e_1,e_2,e_3\}$ and transition probabilities matrix

$$
P=\left(
    \begin{array}{ccc}
      0.5 & 0   & 0.5 \\
      0.3 & 0   & 0.7 \\
      0   & 0   & 1 \\
    \end{array}
  \right)
$$

Then, its associate Markov evolution algebra has the generators $\{e_1,e_2,e_3\},$  the laws
$$
\begin{array}{lllr}
e_1 \cdot e_1 = & 0.5 \, e_1 & + & 0.5 \, e_3, \\
e_2 \cdot e_2 = & 0.3 \, e_1 & + & 0.7 \, e_3, \\
e_3 \cdot e_3 = &            &   &        e_3. \\
\end{array}
$$
\noindent and its weighted directed graph associated is shown in Figure \ref{image2}.

\begin{figure}[ht!]
\centering
\includegraphics[scale=1]{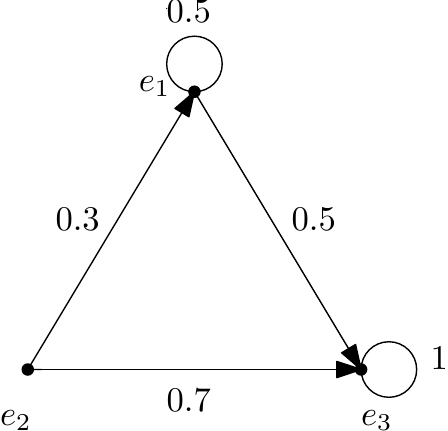}

\caption{The Markov directed graph of Example
\ref{2.5}}\label{image2}
\end{figure}

Obviously, the representation matrix of the evolution operator $L$ is $A=P$. Therefore, $A^{n}=P^{n}$, being $A^{n}=(n_{ij})$ and
$P^{n}=(p_{ij}^{(n)})$, where $p_{ij}^{(n)} = p(X_n = j \mid X_0= i)$ are the $n$-steps transition probabilities.

\vspace{0.1cm}

Particularly, in this example we have

$$
A^2=P^2=
  \left(
    \begin{array}{ccc}
      0.25 & 0   & 0.75 \\
      0.15 & 0   & 0.85 \\
      0    & 0   & 1 \\
    \end{array}
  \right)
    $$

It can be noted that $n_{11}=0.25$ represents the product of the probabilities 0.5 times 0.5, which coincides with the probability that
Markov chain returns state $e_1$ two steps after leaving it, due to the existence of a unique walk from $e_1$ to $e_1$, which is the loop $e_1e_1$.
}
\end{ejem}

\vspace{0.1cm}

We now introduce the following novel concept

\begin{defi} \label{de:weight}
Given a weighted directed graph $G$, the {\em Markov weight} of a walk in $G$ is the product of the weights of the edges of the walk.
\end{defi}

\vspace{0.1cm}

This definition allows us to  state the following theorem

\begin{thm}
Let $G$ be a directed graph with $V(G)=\{x_1,x_2,...,x_n\}$, and let $L$ be the evolution operator of the evolution algebra $A(G)$. If we express
$L^n(e_i)=n_{i1}e_1+n_{i2}e_2+...+n_{ir}e_r$, then $n_{ij}$ coincides with the sum of the Markov weights of the walks of length $n$ from the vertex corresponding to the generator $e_i$ to the vertex corresponding to the generator  $e_j$.
\end{thm}

\vspace{0.1cm}

Previous definition and theorem make possible to improve some
results by Tian. Indeed, Tian showed a result in \cite{book} which
links the concepts of closed subset of states and evolution
subalgebras. It is the following

\vspace{0.1cm}

\begin{thm} [Tian, Theorem 17, Chapter 4 of \cite{book}] \label{cerrsubalg}
Let $C$ be a closed subset of the set of states $S=\{e_i\mid i\in \Lambda\}$ of a Markov chain $X$. $C$ is closed in the sense of probability if and only if $C$ generates an evolution subalgebra of the evolution algebra $M_X$.
\end{thm}

\vspace{0.1cm}

Then, starting from the definition of closed subset of states of a
Markov chain (see Preliminaries), this result can be formulated in
an alternative way, as follows

\vspace{0.1cm}

\begin{thm} \label{th:new1}
Let $C$ be a  subset of the set of generators $\{e_i| i\in \Lambda\}$ of a Markov evolution algebra. $C$ generates an evolution subalgebra of the Markov evolution algebra if and only if
$$n_{ij}=0, \text{ for } e_i \in C, e_j \notin C, n\geq1$$
where $n_{ij}$ represents the elements of the evolution operator $L^n$.
\end{thm}

\vspace{0.1cm}

Tian gave a proof in \cite{book} of the following result (which
can be also proved by 'reductio ad absurdum', taking into
consideration that every evolution subalgebra is an evolution
ideal):

\begin{thm}[Tian, Theorem 18, Chapter 4 of \cite{book}] \label{alsimples}
A Markov chain $X$ is irreducible if and only if its corresponding
evolution algebra $M_X$ is simple.
\end{thm}

Now, as a consequence, we can reformulate that assert as follows

\begin{thm} \label{th:new2}
A Markov evolution algebra is simple if and only if there
exists an integer $n$ verifying that $n_{ij}$, the element of the
representation matrix of the operator $L^n$, is positive for each
pair of generators $e_i$ and $e_j$.
\end{thm}

\vspace{0.1cm}

\begin{nota}
This result allows us to give an alternative way to decide whether
a Markov evolution algebra is simple, without the need of
finding proper ideals of the algebra. Let us see it in the
following example
\end{nota}

\vspace{0.1cm}

\begin{ejem}{\rm

Let $M_X$ be the Markov evolution algebra with generators

 $\{e_1,e_2,e_3,e_4,e_5,e_6\}$ and laws
$$
\begin{array}{lll}
e^2_1 = & 0.3 \, e_2 & + \, 0.7 \, e_6, \\
e^2_2 = & e_3, \\
e^2_3 = & 0.8 \, e_1 &  + \, 0.2 \, e_4, \\
e^2_4 = & e_6, \\
e^2_5 = & e_4, \\
e^2_6 = & e_5. \\
\end{array}
$$

Its associated Markov chain has states $\{e_1,e_2,e_3,e_4,e_5,e_6\}$ and the transition probabilities matrix is

$$
P=\left(
    \begin{array}{cccccc}
      0   & 0.3 & 0 & 0   & 0 & 0.7 \\
      0   & 0   & 1 & 0   & 0 & 0 \\
      0.8 & 0   & 0 & 0.2 & 0 & 0\\
      0   & 0   & 0 & 0   & 0 & 1 \\
      0   & 0   & 0 & 1   & 0 & 0 \\
      0   & 0   & 0 & 0   & 1 & 0 \\
    \end{array}
  \right)
$$

Observe that $P = L,$ $n_{ij}= 0,$ when $i = 3,4,5$ and $j = 1,2,3$ in L and thus for any power of $L.$

Then, according to Theorem \ref{th:new1}, $\{e_4,e_5,e_6\}$ generates an evolution subalgebra (equivalently, it is a proper closed subset in the sense of probability). Therefore, Theorem \ref{th:new2} implies that $M_X$ is not a simple evolution algebra and its associate Markov chain is not irreducible.
}
\end{ejem}

\vspace{0.1cm}

\begin{nota}
Not every subgraph of a weighted directed graph associated with a
Markov evolution algebra generates an evolution subalgebra of
the evolution algebra. For it, the subgraph must meet the
following requirement, which is to be {\em closed}, a new concept
that we now introduce, accompanied by a subsequent example
\end{nota}

\vspace{0.1cm}

\begin{defi} Let $G=(V,E)$ be a
directed graph, not necessarily weighted. We say that the subgraph
$\langle V' \rangle$ induced by $V'\subsetneq V$ is {\em closed}
if $(e_i,e_j)\notin E,$ for all $e_i \in V'$ and $e_j \in
V\setminus V'$.
\end{defi}

\begin{ejem} \label{ejemplo215}
{\rm Let us  consider the directed graph of Figure \ref{image3} (left).

\begin{figure}[ht!]
\centering
\includegraphics[scale=1]{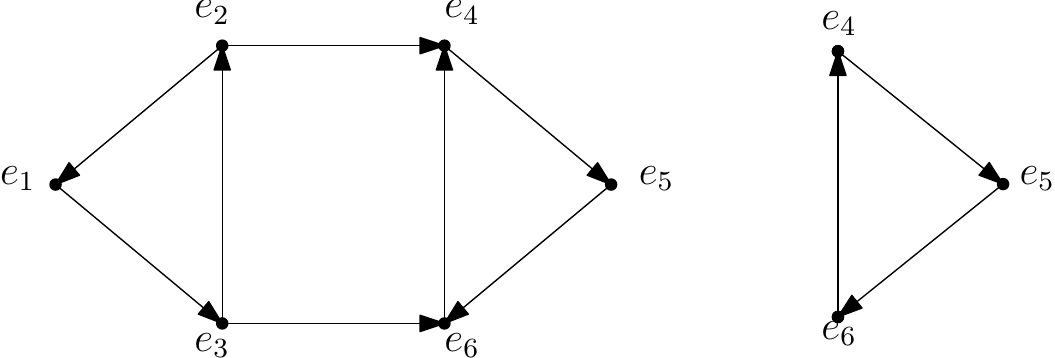}
\caption{The directed graph of Example \ref{ejemplo215}}
\label{image3}
\end{figure}

Recall that the subgraph induced by $\{e_4,e_5,e_6\}$ contains these vertices and those edges between them in the initial graph.  This subgraph closed, since there are no edges $(e_i,e_j),$ with $i=4,5,6$ and $j=1,2,3$. It is shown in Figure \ref{image3} (right).
}

\end{ejem}

This new concept allows us to deduce the following

\begin{thm}\label{th:new3}
Let $G'$ be the weighted directed graph associated to a Markov evolution algebra
$M'$ and let $G$ be the weighted directed graph associated to another Markov evolution algebra
$M$. Then, $G'$ is a closed induced subgraph of $G$ if and only if $M'$ generates an
evolution subalgebra of the evolution algebra $M$.
\end{thm}

\vspace{0.1cm}

This concept also allows us to give a new characteristic of simple
Markov evolution algebras expressed by means of its associated
graph.

\begin{thm}\label{th:new4}
Let $G=(V,E)$ be the directed graph associated to a Markov evolution algebra $M$. Then, $G$ has no closed induced subgraph if and only if the evolution algebra $M$ is simple.
\end{thm}

\vspace{0.1cm}

\subsection{Classification of the generators of a Markov evolution algebra}

Observe that so far there have been several ways to see if a
Markov evolution algebra is simple or not. Now, given a
Markov evolution algebra, we wish to find its simple evolution
subalgebras. To do this, we turn to the classification of
individual generators of a Markov evolution algebra.

\vspace{0.1cm}

Tian
proved in \cite {book} that the concepts algebraically recurrent
generator and recurrent state in the sense of probability are
equivalent, and that the same occurs with transient generator and
transient state in that sense. In addition, Tian proved in Chapter
4 of \cite{book} the following result

\begin{thm} [Tian, Lemma 10 in Chapter 4 of \cite{book}]
A generator $e_i$ is algebraically persistent if and only if all
generators $e_j$ which occurs in $<e_i>$, $e_i$ also occurs in
$<e_j>.$
\end{thm}

This result can be increased and reformulated according to our
research as follows

\begin{thm}  \label{tranrec}
 A generator $e_i$ of the evolution algebra $M$ is transient if and only if there exists one $e_j$ appearing in $\langle e_i \rangle$,
 such that $e_i$ does not appear in $\langle e_j \rangle$. A generator $e_i$ is algebraically recurrent in $M$ if and only if for each generator $e_j$ which appears in $\langle e_i \rangle$, $e_i$ also appears in $\langle e_j \rangle$.
\end{thm}

\vspace{0.1cm}

In this respect, we show next an example in which we classify the
generators of an evolution algebra given.

\vspace{0.1cm}

\begin{ejem}{\rm
Let us see which generators are transient and which recurrent in the Markov evolution algebra defined by

$$
\begin{array}{lllll}
e^2_1 = & 0.4 \, e_1 & + \,  0.2  \, e_2 & + \, 0.2 \,  e_3 & + \, 0.2 \, e_7, \\
e^2_2 = & 0.7 \, e_2 & + \,  0.3  \, e_5, & & \\
e^2_3 = & 0.3 \, e_1 & + \,  0.3  \, e_3 & + \, 0.2 \,  e_4 & + \, 0.1 \, e_5 + \, 0.1 \, e_7,\\
e^2_4 = & 0.9 \, e_4 & + \,  0.1  \, e_7, & & \\
e^2_5 = & 0.5 \, e_2 & + \,  0.5  \, e_5, & & \\
e^2_6 = & 0.3 \, e_2 & + \,  0.5  \, e_4 & + \, 0.1 \,  e_5 & + \, 0.1 \, e_6, \\
e^2_7 = & 0.1 \, e_4 & + \,  0.9  \, e_7. & & \\
\end{array}
$$

To do this, we use Theorem \ref{tranrec} and we obtain that generators $\{e_1, e_3, e_6\}$ are  transient. Analogouly,  generators $\{e_2, e_4, e_5, e_7\}$ are  recurrent.
}

\end{ejem}

Now, we move on to introduce the concept of {\em period} of a
generator $e_i$ of the evolution algebra $M$ in the same way as it
is defined the period of a state $e_i$ in a Markov chain $X$.

\vspace{0.1cm}

\begin{defi}
The {\em period} of a generator $e_j$ of an evolution algebra $M$ is defined as the greatest common divisor of the set of integers $n$ for which $n_{jj}> 0$, where $n_{jj}$ is the $j-th$ diagonal element of the representation matrix of the operator $L^n$. A generator with period 1 is called {\em aperiodic}.
\end{defi}

\vspace{0.1cm}

As a consequence, if $G$ is the weighted directed graph associated with a Markov evolution algebra, the period of a generator $e_j$ is the greatest common divisor of the lengths of the walks from $e_j$ to $e_j$.

\begin{ejem}\label{periodo2}
    {\rm
Given the following Markov evolution algebra
$$e^2_1 = \,e_2, \quad e^2_2 = 0.17 \, e_1 + \, 0.83 \, e_3, \quad e^2_3 = \, e_2,$$

\noindent the weighted directed graph is shown in Figure \ref{image5}.

\begin{figure}[ht!]
\centering
\includegraphics[scale=0.9]{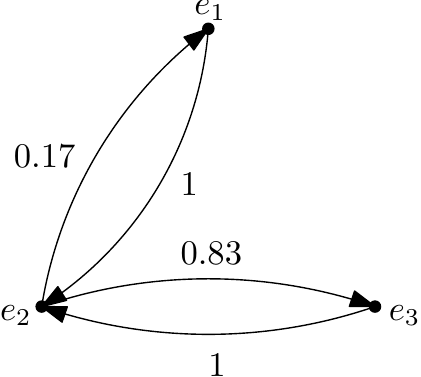}
\caption{The Markov graph of Example \ref{periodo2}}
\label{image5}
\end{figure}

\vspace{0.1cm}

Now, observe that the generator $e_1$ has period 2 because the simple walks $e_1 e_2 e_1$ and $e_1 e_2 e_3 e_2 e_1$ from $e_1$ to $e_1$ have lengths 2 and 4, respectively. Similarly, the period of $e_2$ is the greatest common divisor $\{2,2\}=2$ and the period of $e_3$ is the greatest common divisor $\{2,4\}=2$.
}
\end{ejem}

Note that in this example all generators have the same period. It
is so because in an irreducible Markov chain all
states are either periodic with the same period or aperiodic.
Then, the following result follows

\vspace{0.1cm}

\begin{thm}\label{th:new5}
In a simple evolution algebra all generators are either periodic with the
same period  or all states are aperiodic.
\end{thm}

\vspace{0.1cm}

For instance, according to Theorem \ref{th:new5}, the evolution
algebra of the example \ref{periodo2} is a simple evolution
algebra because its associated Markov chain is irreducible.

\vspace{0.1cm}

In previous results we have dealt with the classification of
individual generators. We will now deal with the classification of
subsets of generators. To do this, we recall that Tian already
proved the following result

\vspace{0.1cm}

\begin{thm} [Tian, Corollary 12, Chapter 4 of \cite{book}]
The state $e_k$ is an absorbent state in a Markov chain $X$ if and
only if $e_k$ is an idempotent element in the evolution algebra
$M_X$.
\end{thm}

Note that this fact is easy to observe since $e_k$ is absorbent in
a Markov chain $X$ if and only if $p_{kk}=1$. So, in the algebra
$M_X$ it is verified $e_k \cdot e_k= e_k$.

\vspace{0.1cm}

Besides, it is easy to note that the absorbent states are closed
subsets of the set of states of a Markov chain and thus, by
Theorem \ref{cerrsubalg}, they generate an evolution subalgebra of
the algebra associated to the chain. It also explains the
following result by Tian

\begin{thm} [Tian, Remark 4 in Chapter 4 of \cite{book}]
Every idempotent element in a Markov evolution algebra $M$
generates an evolution subalgebra of $M.$
\end{thm}

Continuing with our novel  contributions in this study, other concepts which can be introduced in evolution algebras
inspired by their analogous in Markov chain are the following

\vspace{0.1cm}

\begin{defi}
Let $M$ be a Markov evolution algebra. The generator $e_j$ is
said to be {\em accessible} from the generator $e_i$ if there
exists an integer $n$ such that $n_{ij}>0$, where $n_{ij}$ is the
element of the row $i$ and column $j$ of the representation matrix
of the operator $L^n$. Two generators $s_i$ and $s_j$ {\em
communicate} if $s_j$ is accessible from $s_i$ and $s_i$ is
accessible from $s_j$. A generator $e_i$ is said to be a {\em
generator of return} if $a_{ii}>0$, where $a_{ii}$ is the element
of the column $i$ and row $i$ of the representation matrix of the
operator $L$. The set of all generators communicated to generator
$e_i$ constitutes a  {\em class} denoted by $C(e_i)$.
\end{defi}

\vspace{0.1cm}

Then, in the same way as in Markov chains, the generators of an
evolution algebra can be partitioned into connected classes.
Classes may or may not be closed, as it occurs in Markov chains.
If the generator $e_i$ is recurrent, its connected class is
closed, that is, by Theorem \ref{cerrsubalg}, the connected class
to which it belongs generates an evolution subalgebra of the
Markov evolution algebra. Only transient generators can belong
to non closed connected classes.

\vspace{0.1cm}

Besides, if the generator $e_i$ is recurrent and $e_j$ is
accesible from the generator $e_i$, then the generator $e_j$ must
be connected to $e_i$ and must also be recurrent. Thus the
 recurrent generators only connect with other recurrent generators, so
the set of recurrent generators must be closed and therefore
they generate an evolution subalgebra of the Markov evolution
algebra.

\vspace{0.1cm}

Then, if the generator $e_i$ is recurrent, the class $C(e_i)$ is
an irreducible closed set and contains only
 recurrent generators. Thus by Theorems
\ref{cerrsubalg} and \ref{alsimples} if $e_i$ is a recurrent
generator, then $C(e_i)$ generates an evolution subalgebra of the
simple Markov evolution algebra.

\vspace{0.1cm}

These new concepts involve the following  new results in evolution
algebras, which are related with  Corollary 9 in Chapter 3 of
\cite{book}. The first of them is

\begin{thm}
If all generators of an  evolution algebra belong to the same
connected class, then the algebra is simple.
\end{thm}

\vspace{0.1cm}

\noindent whereas the second one is

\vspace{0.1cm}

\begin{thm}\label{partition} The set of
generators of an evolution algebra can be
 partitioned into two subsets. The first one contains only
 transient generators and the second subset only recurrent generators. Moreover,
these last ones generate an evolution subalgebra of the Markov
 evolution algebra and this second subset may also be
 partitioned into connected and irreducible classes which generate simple evolution subalgebras.
\end{thm}

\vspace{0.1cm}

Let us explain these results with an example

\begin{ejem}
\label{ultimo} {\rm Let us consider the $8$-dimensional Markov
evolution algebra
 whose transition probabilities matrix (which coincides with its  evolution operator matrix) is

$$
\left(
    \begin{array}{cccccccccc}
      0.4 & 0.2 & 0   & 0   & 0 & 0.2 & 0.2 & 0 \\
      0.3 & 0.3 & 0.1 & 0   & 0 & 0   & 0.1 & 0.2 \\
      0   & 0   & 0.5 & 0   & 0 & 0.5 & 0   & 0  \\
      0   & 0   & 0.1 & 0.1 & 0 & 0.3 & 0   & 0.5 \\
      0   & 0   & 0   & 0   & 1 & 0   & 0   & 0\\
      0   & 0   & 0.3 & 0   & 0 & 0.7 & 0   & 0\\
      0   & 0   & 0   & 0   & 0 & 0   & 0.8 & 0.2\\
      0   & 0   & 0   & 0   & 0 & 0   & 0.2 & 0.8\\
    \end{array}
  \right)
$$

\vspace{0.15cm}

\noindent and its Markov graph associated $G$ is shown in Figure \ref{image6}.

\begin{figure}[ht!]
\centering
\includegraphics[scale=0.9]{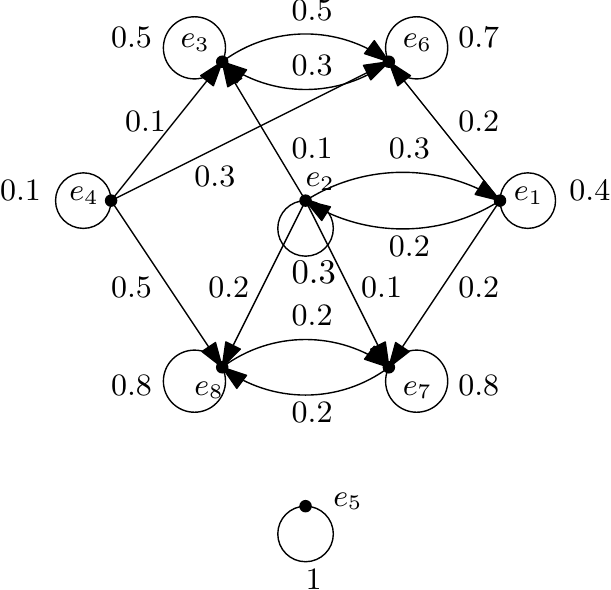}
\caption{The Markov graph of Example \ref{ultimo}}
\label{image6}
\end{figure}

It is easy to observe that the directed graph $G'$ induced by
$\{e_3,e_6\}$ is closed because there are no edges $(e_i,e_j),$
with $i=3,6$ and $j=1,2,4,5,7,8$ in $G$. Therefore, by Corollary
\ref{th:new3} the evolution algebra associated to $G '$ is an
 evolution subalgebra of the initial algebra. In this subalgebra, denoted by
 $E_1$, the generator set is $\{e_3, e_6 \}$ and the law of the algebra is given by

$$
\begin{array}{lll}
e^2_3 = & 0.5 \, e_3 & + \,  0.5  \, e_6\\
e^2_6 = & 0.3 \, e_3 & + \,  0.7  \, e_6\\
\end{array}
$$

Moreover, note that the graph associated to this subalgebra $E_1$
contains no closed induced subgraph. Therefore, by Theorem
\ref{th:new4}, it is a simple evolution subalgebra.

\vspace{0.1cm}

As a conclusion, since every evolution subalgebra is an evolution
ideal of the algebra, the Markov evolution algebra $M$ is not
simple.

\vspace{0.1cm}

 We deal now with the classification of individual generators and
   subsets of the  generator set.

\vspace{0.1cm}

According to Theorem \ref{tranrec} we have that the generators $\{e_1, e_2, e_4\}$ are transient and
generators $\{e_3, e_5, e_6, e_7\}$ are recurrent.

\vspace{0.1cm}

Note that the generators of the evolution algebra $M$ are
 generators of return because the element $a_ {ii}$ of the
 matrix $A$ is positive for each $1 \leq i \leq 8.$ Then, by Theorem
 \ref{partition}, the generator set of $M$ is partitioned into two subsets. The first
 contains only the transient generators $e_1, e_2, e_4$, which form a non closed connected class, and
 the second subset contains the recurrent generators $\{e_3, e_5, e_6,
 e_7, e_8 \}$, which is a closed set. This second subset can be
 partitioned into irreducible and closed connected classes as follows

\begin{itemize}
\item The generator $e_5$ is an idempotent generator. Therefore, it forms a closed irreducible connected class constituted by a unique generator.
\item The connected class of the generator $e_3$ is $\{e_3,e_6\}$, which coincides with the connected class of $e_6$.
This class is closed and irreducible, since it is the connected class
of a recurrent generator. A similar argument works for  generators $\{e_7,e_8\}$.

\end{itemize}

\vspace{0.1cm}

So, the sets $\{e_5\}$, $\{e_3,e_6\}$ y $\{e_7,e_8\}$ generate
three simple evolution subalgebras of $M$.

}
\end{ejem}

\section{Conclusions}

The paper deals with, in the opinion of the authors, the interesting connection between the
Markov chains, the evolution algebras and Graph Theory, that was established
by J. P. Tian when the basic theory of the evolution algebras was introduced
and developed in \cite{book}. Indeed, an homogeneous Markov chain can be canonically
regarded as an evolution algebra. Moreover, an homogeneous
Markov chain (and hence its corresponding evolution algebra), also can be regarded as the weighted directed graph given by the Markov graph associated
to the Markov chain. In this way, the paper is aimed to review some notions by establishing their meaning
in the framework of the evolution algebras, the Markov chains and the graph
theory, providing many examples.

\vspace{0.15cm}

Therefore, the authors think that this research, in addition to being interesting for linking
three areas of Mathematics seemingly disconnected from each other, involves a novel and profound advance in the study of Markov chains.
The advance is due to that this research endows an statistical structure with other two types of different structures, algebraic and discrete,
from which the first one initially lacked and, moreover, introduces new results on them, as all those indicated in the previous section,
which allow us to consider evolution algebras, graphs and Markov chains as a unique mathematical structure.
In this way, it will be indistinct to speak of evolution algebras, graphs and Markov chains, because each property of one of them will immediately imply the same property in the other two.

\vspace{0.15cm}



Our intention in future work is to continue doing research on this topic to step forwards in obtaining new deeper and stronger results that can be interpreted and translated univocally in the frame of these three structures.

\vspace{0.6cm}

\noindent {\em Aknowledgement}: This paper is partially supported
by national projects MTM2015-65397-P, FEDER and autonomic P.A.I.
Groups FQM-326 and FQM-164.


\begin{thebibliography}{99}

\bibitem{Bansaye} Bansaye, V., Ancestral lineages and limit theorems for branching Markov chains in varying environment, {\it J. Theoret. Probab.}
{\bf 32}:1 (2019), 249-281.

\bibitem{bollo} Bollobás, B., Graphs Theory, Springer Verlag, New York, 1979.

\bibitem{cabrera} Cabrera Y., Siles, M. and Velasco, M. V.,  Evolution algebras of arbitrary dimension and their
decompositions, {\em Linear Algebra and its Applications}, {\bf 495}  (2016), 122–162.

\bibitem{csv} Cabrera, Y., Siles, M. and Velasco, M. V., Classification of three-dimensional evolution
algebras, {\it Linear Algebra Appl.} {\bf 524} (2017), 68-108.

\bibitem{cadavid} Cadavid, P., Rodi\~no, M. L. and Rodr\'iguez, P. M., On the connection between evolution algebras,
random walks and graphs. 
{\em  Journal of Algebra and Its Applications}. Vol. 19, No. 02, 2050023 (2020), 28 pages.
 
\bibitem{cadavid2} Cadavid, P., Rodi\~no, M. L. and Rodr\'iguez, P. M., Characterization theorems for the spaces of derivations
of evolution algebras associated to graphs, {\em Linear and Multilinear algebra}. Online 2019: https://doi.org/10.1080/03081087.2018.1541962

\bibitem{cadavid3} Cadavid, P., Rodi\~no, M. L. and Rodr\'iguez, P. M., On the isomorphisms between evolution algebras of graphs and random walks, {\em Linear and Multilinear algebra}. Online 2019: https://doi.org/10.1080/03081087.2019.1645807

\bibitem{elduque} Elduque, A. and Labra, A.,  Evolution algebras and graphs, {\em Journal of Algebra and Its Applications}
{\bf 14}: 7 (2015), 1550103  (10 pages).


\bibitem{harary} Harary, F., Graph Theory, Addison Wesley, Reading, Mass., 1969.

\bibitem{Feller} Feller, W., An introduction to probability theory and its applications. Volume II, Second edition.
John Wiley \& Sons, Inc., New York-London-Sydney, 1971.

\bibitem{GrimSti} Grimmett, G. R. and Stirzaker, D. R., Probability and random processes. Oxford University Press, New York, 2001.

\bibitem{Karlin} Karlin, S. and , Taylor, H. M. An Introduction to Stochastic Modeling, 3rd ed., Academic Press, 1998.

\bibitem{NRV}  Núñez, J., Rodríguez-Arévalo, M. L. and Villar, M. T., Certain particular families of graphicable algebras, {\em
Applied Mathematics and Computation}, {\bf 246}:1 (2014), 416-425.

\bibitem{AMC} Núñez, J., Silvero, M. and Villar, M. T.,
Mathematical tools for the future: Graph Theory and graphicable algebras, {\em Applied Mathematics and Computation},
{\bf 219}:11 (2013) 6113-6125.

\bibitem{Ross} Ross, S. M., Introduction to Probability Models, 10th Ed., Academic Press, 2010.

\bibitem{Staples} Staples, G. S., Graph-theoretic approach to stochastic
integrals with Clifford algebras, {\it J. Theoret. Probab.} {\bf
20}:2 (2007), 257-274.


\bibitem{book} Tian, J. P., Evolution Algebras and their Applications, {\it Lecture Notes in Mathematics}, {\bf Vol 1921} (2008),
Springer-Verlag, Berl\'{\i}n.

\bibitem{invita} Tian, J. P., Invitation to research of new mathematics from biology: evolution algebras.
Topics in functional analysis and algebra, {\it Contemp. Math.,}
{\bf 672} (2016), 257-272.

\bibitem{Petr} Tian, J. P. and Vojtechovsky, P., Mathematical concepts of evolution algebras in non-mendelian genetics,
{\it Quasigroups Related Systems}, {\bf 14}:1 (2006), 111-122.

\end{thebibliography}
\end{document}